# Pre/Post-Assessment Cycles in Calculus: Supporting Student Preparation, Reflection, and Learning

Chamila Malagoda Gamage

Department of Mathematics, University of Florida

cgamage@ufl.edu

**Abstract**

This article describes a semester-long implementation of pre/post-assessment cycles in an undergraduate Calculus I course. The activity was designed to address a common difficulty in lower-division mathematics courses: students often enter a new unit without knowing which prior ideas are relevant, which concepts deserve special attention, or how to judge their own progress before high-stakes exams. Brief pre-assessments were administered before major course units, and related post-assessment questions revisited key ideas after instruction. Students also completed short reflection prompts about what they learned, what remained challenging, and how their understanding changed. The paper analyzes student performance data, participation patterns, and reflective responses to examine how these cycles supported preparation, reflection, and learning. The study contributes a practical model for embedding low-stakes assessment and reflection into the regular structure of an undergraduate mathematics course.

## 1. Introduction

Calculus I is often one of the first university mathematics courses in which students are expected to connect procedures with conceptual reasoning. Students may be able to apply familiar rules for limits, derivatives, and applications, but still have difficulty identifying which prior ideas matter at the beginning of a new unit or how much preparation the new material will require. This difficulty is not only about content knowledge. It also involves preparation, feedback, and self-monitoring. Without a clear preview of the main goals of a new unit, students may not recognize what they do not understand until a homework deadline, quiz, or exam.

This study builds on an earlier implementation of a pre/post-assessment activity in an introductory real analysis course, where a single pre-assessment and post-assessment pair was used to examine student learning around one unit (Gamage, 2025). That earlier work suggested that the activity had value, but it was limited in scope because the process was not repeated throughout the course. The present study is designed as a more structured and sustained version of that idea. Instead of treating pre/post assessment as a one-time activity, the Calculus I implementation uses repeated cycles across the semester so that students encounter assessment, feedback, and reflection as a regular part of the course rhythm. This design allows the paper to examine not only whether students improve after instruction, but also how repeated low-stakes assessment routines may support preparation, self-monitoring, and engagement over time.

## 2. Research Questions

This study is guided by the following research questions:

I. How can pre/post-assessment cycles be integrated into the regular structure of a Calculus I course?

II. How do students describe the role of these cycles in supporting their preparation, reflection, and learning?

III. How can information from pre/post-assessments, and student reflections help instructors identify student needs and adjust instruction?

## 3. Background and Related Literature

### 3.1 Formative Assessment and Feedback

This study is grounded in a formative view of assessment. In this view, assessment is not only a way to measure student performance after instruction, but also a way to support learning while instruction is still taking place. Black and Wiliam (1998) emphasized that formative assessment can provide useful information to both instructors and students when it is used to improve teaching and learning. Similarly, Sadler (1989) argued that feedback is most valuable when it helps students understand the gap between their current performance and the desired goal. Angelo and Cross (1993) similarly emphasize the value of classroom assessment techniques as low-stakes tools that help instructors understand student learning while there is still time to respond. This idea is closely aligned with the present study, where assessment is used not only to record performance but also to guide preparation, feedback, and instructional decisions during the semester.

This perspective is especially relevant in Calculus I, where students often need to connect prior knowledge, new concepts, and critical thinking across several units. A pre/post-assessment cycle can make this process more visible. Pre-assessments can help students notice what they already know and what they need to review, while post-assessments can help them revisit important ideas after instruction. When paired with reflection, these assessments can also support self-monitoring, which Nicol and Macfarlane-Dick (2006) identify as an important part of effective feedback and self-regulated learning.

### 3.2 Knowledge Surveys and Diagnostic Preview

The design of this activity is also related to work on knowledge surveys and diagnostic assessment. Knowledge surveys ask students to judge their confidence or readiness on course-related items, giving both students and instructors a broad view of the course content and learning goals (Nuhfer & Knipp, 2003). Such surveys can help students see the structure of a course and become more aware of what they are expected to learn.

The activity described in this paper differs from a traditional knowledge survey because students answer mathematical questions rather than only rating their confidence. However, it shares an important purpose with knowledge surveys: helping students preview upcoming material and monitor their own learning. In a Calculus I course, this kind of diagnostic preview can be useful because students may enter a new unit with uneven recall of prerequisite ideas. A low-stakes pre-assessment can help make those prerequisite ideas visible before formal instruction begins.

### 3.3 Assessment and Reflection in Undergraduate Mathematics

Recent work in undergraduate mathematics education also supports assessment practices that go beyond checking final answers. Karaali (2018) describes formative assessment as a way to move assessment away from simply ranking students and toward supporting their development as learners. This is closely related to the purpose of the present study, where assessment is intended to help students prepare, reflect, and monitor progress across the semester.

Bermudez and Graysay (2025), use structured pre-assessment as part of a larger instructional design to support students in exploring, developing, and formalizing mathematical arguments. Although their context differs from Calculus I, their work is relevant because it shows how preassessment can be embedded in a planned learning sequence rather than treated as an isolated testing event.

In the calculus context, Jaafar and Lin (2017) describe assessment for learning as a proactive approach to engaging students in the calculus classroom. Their work is relevant to the present study because it supports the use of assessment as part of the learning process, rather than only as a tool for evaluating students after instruction. Reed, Tallman, and Oehrtman (2023) further argue that calculus assessments should attend to the meanings students bring to mathematical ideas, not only to whether students obtain correct final answers. Their work supports the need for assessment items that reveal student thinking and can be used formatively during a course. The present study extends this idea by examining a semester-long structure in which pre-assessments, post-assessment questions, and reflections are used together to support preparation, feedback, and learning in Calculus.

## 4. Course Context and Procedure

### 4.1 Course Context

The activity was implemented in an honors section of Calculus I at a large public university in the United States. Although Calculus I is a large multi-section course, with approximately 1,200 students typically enrolled across all sections in a semester, this study focused on one honors section with 43 students.

The course met three times per week for 50-minute in-person lectures. In addition to lectures, students attended weekly discussion sessions led by teaching assistants. These discussion sessions included an in-person quiz, which was part of the regular course structure. The course also included online homework assignments, midterm exams, and a cumulative final exam.

This setting provided a useful context for implementing pre/post-assessment cycles because the course already included several regular points of student engagement, including lectures, discussion quizzes, homework, and exams. The assessment activity was designed to fit within this existing structure rather than function as a separate or additional major course requirement.

### 4.2 Design of the Pre/Post Assessment Cycles

The pre/post-assessment activity was designed as a semester-long formative routine. Students were informed about the activity, its purpose, and its intended benefits at the beginning of the course. The goal was to help students preview upcoming topics, identify relevant prior knowledge, receive feedback, and reflect on their progress throughout the semester.

There were six online pre-assessment quizzes administered through Canvas, one before each major course unit. The pre-assessments consisted of multiple-choice questions related to ideas that would appear in the upcoming unit. Students were required to complete each pre-assessment before the corresponding unit began. These quizzes were graded for completion only, not for correctness, and together counted for 5% of the final course grade. Although students received credit for submission, the instructor collected the actual scores for comparison with later performance. Upon submission, students received detailed solutions, allowing them to review their responses independently, reflect on their understanding, and prepare for the upcoming lesson if they chose to do so.

After each unit was taught, students completed related post-assessment questions through the regular in-person quiz structure in their discussion sections. These post-assessment questions focused on material students had learned during the unit and were intended to revisit ideas introduced in the pre-assessment. This structure allowed students to encounter related concepts before and after instruction while keeping the activity connected to the normal organization of the course.

The activity also provided useful information for the instructor. Pre-assessment results gave an early view of students' background knowledge and possible misconceptions before instruction began. This allowed the instructor to adjust emphasis, provide additional review when needed, and tailor examples or explanations to student needs. Student progress and perceptions were also gathered through verbal feedback and survey responses, which helped the instructor understand how students experienced the assessment cycles.

### 4.3 How the Activity Fit into the Course

The activity was designed to fit within the existing structure of Calculus I rather than introduce a major new requirement. In large calculus courses, students sometimes complete pre-class activities such as lecture videos or short quizzes before class. Although the pre-assessments in this study were different from lecture quizzes, they used a familiar online format and did not require a major change to the course organization.

The post-assessment component also fits naturally into the course because all sections already included weekly in-person discussion quizzes. Related post-assessment questions were incorporated into this existing quiz structure. As a result, students encountered the activity as part of the normal rhythm of the course: previewing ideas before a unit, learning the material in lecture and discussion, revisiting related questions after instruction, and reflecting on their progress.

This design kept the activity at low stakes and manageable. The pre-assessments were completed online and graded for submission, while the post-assessment questions were connected to assessments students were already taking. This allowed the instructor to gather useful information about student readiness and progress without substantially increasing class time, grading demands, or disruption to the course schedule.

### 4.4 Data Sources

The study used both quantitative and qualitative data from the course. The quantitative data included student completion records and scores from the six Canvas pre-assessments, as well as scores on related post-assessment questions from the weekly discussion quizzes. Although the pre-assessments were graded for completion only, the actual scores were recorded for comparison.

The qualitative data included student reflection responses, survey feedback, and verbal feedback collected during the semester. These responses provided information about how students experienced the activity, what they found useful, and what remained challenging.

## 5. Data Analysis

### 5.1 Quantitative Results

The quantitative analysis used gradebook data from 43 students in the honors Calculus I section. There were six pre-assessments and twelve post-assessment quizzes, and since they had different point values, all scores were converted to percentages before analysis.

Each pre-assessment was paired with two related post-assessment quizzes. For each unit, the post-assessment score was calculated as the average of the two related quiz percentages. The overall pre-assessment score was calculated as each student's average across the six pre-assessments, and the overall post-assessment score was calculated as each student's average across the twelve post-assessment quizzes.

Participation in the assessment cycles was very strong. All 43 students completed all six pre-assessments and all twelve post-assessment quizzes. Students also earned relatively high scores on the pre-assessments, with an overall mean of 69.86% and median of 71.85%. Since these quizzes were graded for completion rather than correctness, this pattern suggests that students generally engaged with the questions in a meaningful way.

To examine changes from pre-assessment to post-assessment, paired-samples t-tests were conducted. Cohen's d for paired samples was also calculated to estimate the size of the change. Correlations were used to examine the relationship between students' average pre-assessment scores, average post-assessment scores, and exam performance.

### 5.1.1 Overall Pre/Post Performance

Across the semester, students' average pre-assessment score was 69.86%, while their average post-assessment score was 88.61%. This represents an average increase of 18.75 percentage points. A paired-samples *t*-test showed that the increase from pre-assessment to post-assessment was statistically significant. This means that the observed increase was unlikely to have occurred by chance. The effect size was large, indicating a substantial improvement from pre- to post-assessment performance.The standard deviation also decreased from 14.85 on the pre-assessments to 8.33 on the post-assessments, suggesting that students' scores were more spread out before instruction, possibly due to uneven prior knowledge, and more consistent after instruction.

| Measure | Pre-Assessment | Post-Assessment |
|---|---|---|
| Number of students | 43 | 43 |
| Mean score | 69.86% | 88.61% |
| Standard deviation | 14.85 | 8.33 |
| Median score | 71.85% | 90.83% |
| Paired *t*-test | 9.45, $p<0.001$ | |
| Cohen's *d* | 1.44 | |

Table 1: Overall performance

### 5.1.2 Unit-by-Unit Performance

The unit-by-unit comparisons show that students improved from pre-assessment to post-assessment in all six units. All six gains were statistically significant. The largest gain occurred in Unit 4, where the pre-assessment average was the lowest. This suggests that the activity was especially useful in identifying a unit where students entered with weaker initial understanding.

| Unit | Pre-Mean | Post-Mean | Paired *t* | Cohen's *d* |
|---|---|---|---|---|
| 1 | 71.51% | 87.62% | 5.86 | 0.89 |
| 2 | 71.51% | 91.05% | 5.29 | 0.81 |
| 3 | 77.33% | 85.47% | 3.03 | 0.46 |
| 4 | 52.56% | 93.14% | 10.14 | 1.55 |
| 5 | 72.09% | 86.98% | 3.98 | 0.61 |
| 6 | 74.16% | 87.44% | 3.30 | 0.50 |

Table 2: Unit-by-unit performance

## 5.2 Qualitative Results

Although the pre-assessments were graded for completion rather than correctness, students appeared to engage with them seriously. The overall pre-assessment average was approximately 70%, suggesting that students did more than simply submit answers for credit. This is important because one possible concern with completion-based quizzes is that students may not put meaningful effort into them. In this course, the relatively high pre-assessment performance suggests that many students used the quizzes as intended: to preview upcoming material and assess their own readiness.

Student feedback suggested that the pre-assessments helped students enter class with greater awareness of upcoming topics. Before a unit began, students had already seen relevant terminology, formulas, and types of questions. This did not mean that they were expected to master the content independently before lecture. Rather, the pre-assessments gave them an initial exposure to the ideas that would later be developed in class.

Student responses gave a more direct view of how students experienced this preparation. Several comments emphasized that the most useful part of the activity was seeing the upcoming material before formal instruction. For example, one student wrote that the activity helped by "*knowing what type of material I would be learning,*" while another described the most useful part as "*seeing what was coming up*." A third student gave a more detailed explanation: "*I liked how I could see topics for the next unit, which helped me know what we would be doing and clarify what I know and what I don't know.*" These responses suggest that, for some students, the pre-assessments helped reduce uncertainty at the beginning of a unit and gave them a clearer sense of what to expect.

Other responses connected the activity to self-monitoring. One student identified the most useful part as "*understanding where my knowledge currently sits*." This comment aligns with the purpose of the activity: to help students recognize their current understanding before instruction and then compare that understanding with their later performance. In this way, the pre/post-assessment cycle supported reflection, not only performance.

The student feedback was not uniformly positive, which is also important to note. Some students indicated that the pre-assessment did not substantially change how they studied or that it was only somewhat helpful. This suggests that students may use pre-assessments in different ways. Some used it as an opportunity to preview the unit and check their understanding, while others seemed to engage with it more lightly.

The repeated nature of the activity was an important part of the student experience. Because the pre/post-assessment cycle occurred throughout the semester, students encountered a regular pattern of previewing material, learning it in class, revisiting related questions, and reflecting on their progress. This routine may have helped students maintain more consistent effort across the course.

This is especially relevant in a semester-long calculus course, where students can lose momentum near the end of the semester. The assessment cycles gave students recurring checkpoints and encouraged them to stay engaged with each unit as it began. The reflection prompts also helped shift the activity from simple performance measurement to self-monitoring. Students were asked to think about what they learned, what remained challenging, and which strategies helped them improve.

The activity also provided useful information for the instructor. Before beginning a unit, the instructor could review pre-assessment results to identify common mistakes, gaps in prior knowledge, or topics that students already understood relatively well.

This information helped the instructor decide where to spend more time, what examples to emphasize, and which prerequisite ideas needed review.

For example, the unit-by-unit quantitative results showed the largest gain in Unit 4, which focused on applications of derivatives. The low pre-assessment performance in this unit suggested that students either had weaker prior knowledge or found the topic more difficult at the start. Having this information before instruction allowed the instructor to give additional attention to the unit and address common difficulties more directly.

## 6. Discussion

The results suggest that the pre/post-assessment cycles functioned as a useful formative structure in the Calculus I course. The main value of the activity was not simply that students scored higher on post-assessments, since improvement after instruction is expected. Rather, the cycles created a regular process in which students previewed upcoming material, received feedback, revisited related ideas, and reflected on their progress. This routine supported preparation before class and encouraged more consistent engagement across the semester.

The activity also provided timely information for the instructor. Pre-assessment results made it possible to identify common mistakes and gaps in prior knowledge before a unit began. This was important because the information could be used while the same students were still learning the material, rather than only after the semester had ended. In many courses, assessment data are reviewed too late to benefit the students who generated them. In this activity, however, the instructor could respond in real time by adjusting emphasis, reviewing prerequisite material when needed, and spending more time on topics where students appeared less prepared. For example, the lower pre-assessment performance and larger gain in Unit 4 suggested that applications of derivatives required additional attention, allowing the instructor to adjust instruction and dedicate more time to this topic during the unit.

This study is based on a previous intervention involving a single pre/post-assessment pair in an introductory real analysis course (Gamage, 2025). That earlier activity suggested potential benefits, but it also raised practical concerns. If a pre-assessment is announced in advance, student preparation may affect its diagnostic value; if it is not counted toward the grade, students may not take it seriously. In addition, an in-class pre-assessment can reduce instructional time and exclude students who are absent. The present Calculus I activity addressed these issues by embedding the pre-assessments into the regular course structure, making them a small completion-based part of the course grade, and administering them online. As a result, all students had access to the pre-assessments, the activity did not take lecture time, and students encountered the process as a normal part of the course rhythm.

Overall, the quantitative and qualitative findings support the value of studying pre/post-assessment cycles as a semester-long formative routine. The results should not be interpreted as proving that the assessment cycles alone caused student learning. A stronger interpretation is that the cycles provided a consistent structure for preparation, feedback, review, and reflection within the existing Calculus I course.

## 7. Limitations and Future Directions

### 7.1 Limitations

This study has several limitations. First, the activity was implemented in one honors section of Calculus I with 43 students. Since this was a relatively small and selected group of students, the findings may not generalize directly to all Calculus I sections, especially larger or non-honors sections. Student readiness for this type of activity may also vary. Some students may welcome the opportunity to preview material before instruction, while others may feel uncomfortable or anxious when asked to engage with unfamiliar ideas on their own.

Second, although other Calculus I sections were offered during the same semester, this study did not use a formal control group. Therefore, the results should not be interpreted as evidence that the pre/post-assessment cycles alone caused the observed gains. Students' improvement was also influenced by regular instruction, homework, discussion quizzes, studying, and other course experiences. The purpose of the study was not to isolate the activity as the only cause of learning, but to examine how it functioned as a formative structure within the course.

Third, scalability remains an important challenge. Calculus I is a large multi-section course, and implementing a semester-long assessment cycle across all sections would require coordination among instructors, teaching assistants, and course administrators. Even when the pre-assessments are online and graded for completion, the design of meaningful questions, alignment with discussion quizzes, and use of student feedback require careful planning.

**7.2 Future Directions**

Future work could extend this activity to larger lecture sections and multiple sections of Calculus I. A larger implementation would make it possible to compare different sections more systematically and to study whether similar patterns appear across different instructors, student populations, and course formats. Such work would require conversation and coordination among instructors so that the activity is implemented consistently while still allowing flexibility for individual teaching styles.

The assessment design could also be refined. Future versions of the pre-assessments could include different types of items, such as conceptual questions, procedural questions, application problems, and questions targeting common misconceptions. In addition, students could be asked to rate their confidence on selected items. This would allow the instructor to compare what students know with how confident they feel, giving a more complete picture of student readiness.

Future research could also examine longer-term outcomes. For example, it would be useful to study whether students who experience repeated pre/post-assessment cycles retain calculus concepts more effectively, feel more prepared for later mathematics courses, or continue to use similar self-monitoring strategies in future studies. These questions would help determine whether the benefits of the activity extend beyond a single semester.

## 8. Conclusion

This study examined a semester-long pre/post-assessment cycle in an honors Calculus I course. The results showed strong student participation, significant improvement from pre-assessment to post-assessment across all units, and a positive relationship between post-assessment performance and exam performance. More importantly, the activity created a regular structure for students to preview upcoming ideas, receive feedback, revisit important concepts, and reflect on their learning.

The findings suggest that pre/post-assessment cycles can be more than a way to measure progress. When embedded into the normal rhythm of a course, they can support preparation, self-monitoring, and instructional responsiveness without adding a major burden to the course structure. While the results do not prove that the assessment cycles alone caused student learning, they show that a carefully designed, low-stakes formative routine can help make learning more visible for both students and instructors.

**Declaration**

Funding: This research received no external funding.

Conflicts of Interest: The author declares that there are no conflicts of interest.

Data Availability: The data sets used and analyzed during the current study are available from the corresponding author upon request.